# Contradictions in some primes conjectures

VICTOR VOLFSON

ABSTRACT. This paper demonstrates that from the Cramer's, Hardy-Littlewood's and Bateman-Horn's conjectures (suggest that the probability of a large positive integer $x$ being a prime - $1/\ln(x)$ ) it follows that the events consisting in a positive integer $x$ being not divisible by different primes are dependent with the ratio $0,5e^{\gamma}$ ($\gamma$ - Euler's constant). In establishing Hardy-Littlewood's and Bateman-Horn's conjectures, their authors followed the first suggestion by another one assuming the independence of the above-mentioned events, which on the basis of the first suggestion and Mertens' theorem is not exactly. This paper demonstrates why these suggestions do not lead to an erroneous result using the Hardy-Littlewood's conjecture for twin primes as an example. The author provides generalized conjectures, which if taken together with the first suggestion, make Hardy-Littlewood's prime $k$ - tuples and Bateman-Horn's conjectures true.

1. INTRODUCTION

Hardy-Littlewood's conjecture [1] Bateman-Horn's conjecture [2] and Cramer's conjecture [3], which are based on probabilistic models suggest that the probability of the event that a large positive integer $x$ is prime, is

$$Pr(A_1) = 1/\ln(x)$$ (1.1).

The first idea occurring after this assumption suggests that if $x$ is known, then the probability of its being a prime may be equal either to 0 or to 1. Perhaps that is why it is assumed that the exact value of $x$ is unknown. The authors of the said conjectures certainly do not consider (1.1) to be correct, but only suggest that since proceeding from the law of primes, their asymptotic density equals (1.1). However, only the density of primes on a finite interval is a probability, whereas an asymptotic density is not [4]. Thus, this assumption in general is not correct. Despite this, let us proceed from this assumption of the authors of the conjectures and see where it takes us. From now on the paper is written for the instance when probabilistic assumption (1.1) holds true, and its aim is to show what its result will be.

Keywords: density of the sequence, probability, contradictions, twin primes conjecture, prime $k$-tuples conjecture.

___________________

## 2. CONTRADICTIONS OF CONJECTURES

Thus, if the probability assumption (1.1) holds, this means:

1. There is a chance of a particular natural number $x$ there be prime.

2. The probability of this chance is numerically equal $1/\ln(x)$.

Based on the law of the prime and the assumption (1.1) the probability of a large natural number $x$ being a prime is the asymptotic density of primes. In this formal sense we understand the unconditional and conditional probability of large natural number being a prime in this paper.

Let us formulate the Hardy-Littlewood's conjecture through the asymptotic density of the sequence for twin primes.

The asymptotic density of the twin prime sequence $f(n)$ of positive integers $x \geq 2$ equals $P(f, 2, x) \sim C/\ln^2(x)$, where $C = 2 \prod_{p>2} \frac{1-2/p}{(1-1/p)^2}$.

Proceeding from the Marten's theorem, the following asymptotic correlation will be true:

$$1/\ln(x) \sim 0,5e^\gamma \prod_{p \leq \sqrt{x}}(1 - 1/p)$$, (2.1) where $\gamma$- Euler's constant.

Let us proceed from the fact that if the positive integer x is a prime number, it is not divisible by the primes $p$, for which $p \leq \sqrt{x}$ will hold. Among $p$ consecutive integers only one divisible by $p$, so the probability of the opposite event that a natural number $x$ is not divisible by p is equal to

$$1 - 1/p \quad (2.2)$$

As the probabilities of the events that the number x is not divisible by different primes, for an infinite number of events, may be dependent, the probability that an unknown positive integer $x$ is the prime $A_1$ on the basis of (2.2) may be written as:

$$Pr(A_1) \sim C_1 \cdot \prod_{p \leq \sqrt{x}}(1 - 1/p)$$, (2.3) $C_1$- dependency ratio.

On the other hand, in view of the probability assumption (1.1) and formula (2.1):

$$Pr(A_1) = 1/\ln(x) \sim 0,5e^\gamma \cdot \prod_{p \leq \sqrt{x}}(1 - 1/p)$$ (2.4)

Comparing (2.3) and (2.4) we will obtain: $C_1 = 0,5e^\gamma$. Thus, the events that unknown large positive integer $x$ is not divisible by different primes, when an infinite number of primes, are dependent with the dependency ratio $C_1 = 0,5e^\gamma$, where $\gamma$ - Euler's constant and $0,5e^\gamma = 0,8905362...$.

However, Hardy-Littlewood's and Bateman-Horn's conjectures on primes suggest that these events are independent, which, as I have demonstrated, is not exactly if the conjecture premise holds (1.1).

The question naturally arises: if such contradiction exists, are the conjectures themselves true? And if they are true – why so? I shall try to answer this question.

Let's consider the probability of the $A_2$, event in which an unknown integer $x+2$ is prime, in view of the probability assumption (1.1). Taking into account the asymptotics, i.e. that $x$ approaches infinity, proceeding from (2.4) we shall get:

$$Pr(A_2) = 1/\ln(x) \sim 0,5e^\gamma \cdot \prod_{p \leq \sqrt{x}}(1 - 1/p)$$  (2.5)

In view of the probability assumption (1.1) and on the basis of formulas (2.4) and (2.5), the Hardy-Littlewood's conjecture for twin primes may be written in a probabilistic form – the probability that the positive integers $x$ and $x+2$ are simultaneously primes will equal to:

$$Pr(A_1 \cdot A_2) = Pr(A_1)Pr(A_2/A_1) = C \cdot Pr(A_1) \cdot Pr(A_2) = C/\ln^2(x),$$

where $Pr(A_2/A_1)$ - probability that the unknown positive integer $x+2$ will be a prime provided $x$ is already a prime, and $C$ - dependency ratio.

Then let us find the dependency ratio - $C = Pr(A_2/A_1)/Pr(A_2)$.

We divide all the natural numbers into classes modulo $p$, where $p$ is an arbitrary prime. There will be $p$ such classes, each of which will have the same amount of numbers. If a natural number $x$ is not divisible by $p$, it must be in one of the $p-1$ classes. Of these, the number of classes, where the integer $x-2$ is not divisible by $p$ is equal to $p-2$. Therefore the probability of the event that the natural numbers $x$ and $x-2$ is not divisible by p is equal to: $(p-2)/(p-1)$ (2.6), except for the case $p=2$, when the probability equals to 1.

In view of (2.6) and the possible dependency for an infinite number of events, that the positive integer $x+2$ is not divisible by different primes, provided that the positive integer $x$ is a prime, the unknown probability will be found by:

$$Pr(A_2/A_1) \sim C_2 \prod_{2<p\leq\sqrt{x}} \frac{p-2}{p-1}.$$  (2.7)

I would suggest that the events consisting in the large positive integer $x+2$ being not divisible by different primes, when an infinite number of primes, are equally dependent whether or not the positive integer $x$ is prime, i.e. the dependency ratio in both cases equals to $0,5e^{\gamma}$ (2.8).

Thus, proceeding from my conjecture (2.8) - $C_2 = 0,5e^{\gamma}$, that is why proceeding from formulas (2.7), (2.8), we will get:

$$Pr(A_2/A_1) \sim 0,5e^{\gamma} \prod_{2<p\leq\sqrt{x}} \frac{p-2}{p-1} = 0,5e^{\gamma} \prod_{2<p\leq\sqrt{x}} \frac{1-2/p}{1-1/p}.$$  (2.9)

Considering formula (2.9), the dependency ratio will make:

$$C(x) = Pr(A_2/A_1)/Pr(A_2) \sim \frac{\prod_{2<p\leq\sqrt{x}}(p-2)/(p-1)}{\prod_{p\leq\sqrt{x}}(p-1)/(p)} = 2\prod_{2<p\leq\sqrt{x}} \frac{p(p-2)}{(p-1)^2}.$$  (2.10)

Thus, the $0,5e^{\gamma}$ ratios in (2.10) decrease and do not affect the ultimate result:

$$C = 2\lim_{x\to\infty} \prod_{2<p\leq\sqrt{x}} \frac{p(p-2)}{(p-1)^2} = 2\prod_{p>2} \frac{p(p-2)}{(p-1)^2} = 2\prod_{p>2} \frac{1-2/p}{(1-1/p)^2},$$  (2.11)

which corresponds to the Hardy-Littlewood;s conjecture for twin primes.

Now let's consider the Hardy-Littlewood's conjecture for prime $k$-tuples. Let's formulate the Hardy-Littlewood's conjecture through the asymptotic density of the sequence for prime $k$-tuples.

The asymptotic density of the sequence of $f_k(n)$ prime $k$-tuples on the $x\geq 2$ interval for the natural sequence will make: $P(f_k,2,x) \sim D_k/\ln^k(x)$, where $D_k = \prod_{p>2} \frac{1-w_k(p)}{(1-1/p)^k}$ and $w_k(p)$ is the number of compared solutions:

$$x(x+2n_1)...(x+2n_1+...+n_{k-1}) \equiv 0(\mod p).$$  (2.12)

I would like to state a more general conjecture that the events consisting in a large positive integer: $x+2n_1+...x+2n_k$ being not divisible into different primes, when an infinite number of primes, are equally dependent whether each of the numbers $x, x + 2n_1, ..., x + 2n_1 + ...2n_{k-2}$ is whether a prime or not, i.e. the dependency ratio in each instance equals to $0,5e^\gamma$ (2.13).

It is possible to show that if premise (1.1) and conjecture (2.13) hold, conjecture (2.12) for prime $k$-tuples will be true.

Another conjecture based on the probabilistic assumption (1.1) is the Bateman-Horn's conjecture.

Let us formulate this conjecture also through the asymptotic density of the sequence of positive integers shown below. Let $g_1, g_2, ...g_k$ be sequences of non-reducible polynomials with integer coefficients and the $h_1, h_2, ...h_k$ power, correspondingly, each of them taking an infinite number of prime values. Then the asymptotic density of the sequence of positive integers on the interval of the $x \geq 2$ natural sequence, at which all polynomials take prime values - $g(n)$ is defined by the:

$$P(g, 2, x) \sim \frac{E_k}{H_k \ln^k(x)}, \quad E_k = \prod_p \frac{1 - \alpha_k(p)/p}{(1 - 1/p)^k}$$ formula (2.14), where $H_k = h_1 \cdot h_2 \cdot ... \cdot h_k$ and $\alpha_k(p)$ is the number of compared solutions:

$$g_1(x) \cdot g_2(x) \cdot ... \cdot g_k(x) \equiv (0 \mod p).$$

Let's generalize conjecture (2.13) over the Bateman-Horn's conjecture. The events that the positive integers: $g_1(x), g_2(x)..g_k(x)$ are not divisible by different primes, when an infinite number of primes, are equally dependent if $g_1(x), g_2(x)..g_k(x)$ are prime numbers or not, i.e. the dependency ratio in all the instances indicated equals $0,5e^\gamma$. (2.15)

It is also possible to show that if the probability assumption (2.1) and conjecture (2.15) hold, the Bateman-Horn's conjecture (2.14) will also hold true.

3. CONCLUSIONS AND SUGGESTIONS FOR FURTHER WORK

This work demonstrates that the probability assumption (1) in the Cramer's, Hardy-Littlewood's and Bateman-Horn's conjectures results in a dependency of the events that the unknown positive integer $x$ is not divisible by different primes (an infinite number of primes). The ratio of this dependency equals to $0,5e^\gamma$, where is $\gamma$ Euler's constant.

When formulating Hardy-Littlewood's and Bateman-Horn's conjectures their authors follow the probability assumption (1.1) by another assumption of independence of the above events, which is not exactly. This is the contradiction in the data output given conjectures. Grenville [5] pointed out the contradiction in the similar Cramer's conjecture.

The work demonstrates why these assumptions do not lead to an erroneous result using the Hardy-Littlewood's conjecture for twin primes as an example.

The author provides conjectures, which, once true, simultaneously with the probability assumption (1.1), will make Hardy-Littlewood's prime $k$ - tuples and Bateman-Horn's conjectures hold.

In the next paper the author intends to demonstrate that Hardy-Littlewood's prime $k$ - tuples and Bateman-Horn's conjectures will hold if the probability assumption (1.1) and the generalized conjectures provided in the work hold true.

4. ACKNOWLEDGEMENTS

I would like to thank all those who participated in the discussion of this work and express my gratitude for their invaluable comments and suggestions for its further improvement.